\documentclass[10pt]{article}
\usepackage{amssymb,amsmath,amsthm,epic,eepic}
\parskip=\medskipamount
\font\fr=eufm10  scaled \magstep 1   %Caracteres g¢ticos.
 1  %Caracteres "doble palo".

\def\cal#1{\mathcal{#1}}

\def\t{\tau}
\def\s{\sigma}
\def\w{\omega}

\def\G{\Gamma}
\def\g{\gamma}

\def\e{\epsilon}

\def\cinfty#1{C^{\scriptscriptstyle\infty}(#1)}
\def\vectorfields#1{{\cal X}(#1)}

\def\fpd#1#2{\frac{\partial #1}{\partial #2}}
\def\R{{\rm I\kern-.20em R}}
\def\F{{\rm I\kern-.20em F}}

\def\ovl#1{\overline{#1}}

\def\tilde#1{\widetilde{#1}}
\def\alg#1{\mbox{\fr #1}}

\def\half{\mbox{$\frac{1}{2}$}}

\DeclareMathOperator{\sign}{sgn} \DeclareMathOperator{\spn}{span}
\DeclareMathOperator{\dom}{dom} \DeclareMathOperator{\im}{im}
\def\rel#1#2{#1 \leftrightarrow #2}
\def\ot{\leftarrow}

\def\oneforms#1{{\cal X}^*(#1)}

\def\ro{\rho}

\def\sign{\mathrm{sgn}}
\def\etal{et al.}

\title{Leafwise holonomy of connections over a bundle map}
\author{B. Langerock \\ Department of Mathematical Physics and Astronomy,\\ Ghent
University, Krijgslaan 281 S9,B-9000 Ghent Belgium \\
{\ttfamily \small bavo.langerock@rug.ac.be}\\
{\ttfamily \small http://maphyast.rug.ac.be}}

\newtheorem{thm}{Theorem}
\newtheorem{cor}[thm]{Corollary}
\newtheorem{lem}[thm]{Lemma}

\newtheorem{prop}[thm]{Proposition}

\theoremstyle{definition}
\newtheorem{defn}[thm]{Definition}
\newtheorem{exmp}[thm]{Example}
\newtheorem{rem}[thm]{Remark}

\def\proofname{Proof. }
\newenvironment{pf}{{\bfseries\proofname}}

\begin{document}
\maketitle

\begin{abstract}
In this paper we introduce a generalisation of the notion of
holonomy for connections over a bundle map on a principal fibre
bundle. We prove that, as in the standard theory on principal
connections, the holonomy groups are Lie subgroups of the
structure group of the principle fibre bundle and we also derive a
straightforward generalisation of the Reduction Theorem.
\end{abstract}

\noindent {\bf Keywords:}  holonomy groups, generalised
connections

\noindent {\bf 2000 MSC:} 53C05, 53C29

\section{Introduction}\label{S1}
The standard notion of a connection (see e.g. \cite{koba}) has
been generalised along many different lines. To mention a few of
these generalisations, we may refer to the so-called
``pseudo-connections'' (see F. Etayo \cite{Etayo} for a review)
and, in particular, to the ``partial connections'' as studied, for
instance, by F. Kamber \etal\ \cite{KamTon}. More recently, a
notion of connection on Lie algebroids has been introduced and
studied, among others, by R.L. Fernandes \cite{Fern1,Fern2}. The
importance of these generalisations can be illustrated, for
instance, by the fact that partial connections were used to prove
the vanishing of some cohomology classes on manifolds admitting a
regular integrable distribution, and the theory of Lie algebroid
connections has lead to the construction of a generalised
Chern-Weil homomorphism onto the set of Lie algebroid cohomology
classes.

In an attempt to establish a unified framework for the various
types of connections mentioned above, we have introduced in a
recent paper a notion of generalised connections over a vector
bundle map \cite{algcon}. In some subsequent papers, we have
investigated possible applications of these generalised
connections: to nonholonomic mechanics \cite{mijzelf}, to the
study of length minimising curves in sub-Riemannian geometry
\cite{sub} and to the formulation and proof of a geometric version
of the maximum principle in control theory \cite{control}.

In \cite{algcon} we managed to associate a notion of
``parallelism'' and ``covariant derivation'' with a generalised
connection over a bundle map. However, torsion and curvature are
in general not well defined unless the bundles under consideration
admit some additional geometric structures, such as in the case of
a pre-Lie algebroid. In this paper we present a notion of
``holonomy'' for these generalised connections and we derive a
version of the Reduction Theorem \cite[p 83]{koba}. It should be
mentioned that holonomy has already been studied for partial
connections in the framework of (contact) sub-Riemannian geometry,
see for instance \cite{falbel}, and for generalised connections in
the framework of Lie algebroids \cite{Fern2}.

The structure of the paper is as follows. In Section
\ref{sectionanchor} we introduce the notion of an anchored bundle
and discuss some of its basic properties. The structure of an
anchored bundle is encountered in sub-Riemannian geometry, control
theory, nonholonomic mechanics and also in the theory of (affine)
Lie algebroids \cite{grab1,tomaffgebroid,tomwil}. Therefore, we
believe it is worth to study this structure in its own right. In
Section \ref{sectieprincipalcon} we introduce the notion of a lift
over an anchored bundle, which can be regarded as a right
invariant anchor map on a principal fibre bundle, commuting with a
given anchor map on the base space. Furthermore, the notion of
``leafwise holonomy'' of a lift over an anchored bundle is
defined. In Section \ref{sectieholonomy} we prove that the
generalised holonomy groups are Lie subgroups of the structure
group of the given principal fibre bundle. A generalisation of the
Reduction Theorem is then easily obtained. To conclude this paper
we briefly discuss some possible applications in sub-Riemannian
geometry.

All manifolds considered in this paper are real, finite
dimensional smooth manifolds without boundary, and by smooth we
will always mean of class $C^{\scriptscriptstyle\infty}$. The set
of (real valued) smooth functions on a manifold $B$ will be
denoted by $\cinfty{B}$, the set of smooth vector fields by
$\vectorfields{B}$ and the set of smooth one-forms by
$\oneforms{B}$. The set of all smooth (local or global) sections
of an arbitrary fibre bundle $\tau: E \rightarrow B$ will be
denoted by $\G(\tau)$.

\section{Anchored bundles} \label{sectionanchor}

In this section we describe the basic structure on which our study
of generalised connections is based, namely that of an anchored
bundle. Let $M$ denote an arbitrary $n$-dimensional manifold with
tangent bundle $\t_M : TM \to M$. The conceptual idea of an
anchored bundle is that one considers a bundle over $M$ which is
related to $TM$, in such a way that, for further developments, the
bundle can be taken as an alternative to the tangent bundle of
$M$. The notion of an anchored bundle is also encountered in the
work of P. Popescu \cite{pop92}, who also uses the denomination
``relative tangent space''.

\begin{defn}
An anchored bundle on $M$ is a pair $(\nu,\ro)$ where, $\nu:N\to
M$ denotes a fibre bundle over $M$, and $\ro:N \to T M$ is a
bundle map, fibred over the identity on $M$. We call $\ro$ the
{\em anchor map} of the anchored bundle.
\end{defn}

The following diagram is commutative:

\setlength{\unitlength}{.7cm}
\begin{center}
\begin{picture}(16,3)(-3,5.5)
\thicklines \put(2.3,8){$N$} \put(6.5,8){$TM$} \put(4.4,5.4){$M$}
\put(3,8.2){\vector(1,0){3.2}} \put(2.7,7.6){\vector(1,-1){1.5}}
\put(6.6,7.6){\vector(-1,-1){1.5}} \put(2.9,6.6){$\nu$}
\put(6.2,6.6){$\tau_M$} \put(4.4,8.4){$\rho$}
\end{picture}
\end{center}

We say that an anchored bundle $(\nu,\ro)$ is {\em linear}, if
$\nu$ is a vector bundle and $\ro$ is a linear bundle morphism.

Consider two anchored bundles $(\nu',\ro')$ and $(\nu,\ro)$ with
base manifolds respectively $M'$ and $M$. An anchored bundle
morphism $(f,\ovl f)$ from $(\nu',\ro')$ to $(\nu,\ro)$ consists
of a smooth mapping $\ovl f: M' \to M$ and a bundle morphism $f:
N' \to N$ fibred over $\ovl f$, in such a way that the following
equation holds:
\[
T \ovl f \circ \ro' = \ro \circ f.
\]
We say that $f$ is an anchored bundle isomorphism if $f$ is a
bundle isomorphism (see e.g. \cite{saunders}), and if, in
addition, $f^{-1}$ is also an anchored bundle morphism. In this
case we can write $\ro' = T (\ovl f)^{-1}\circ \ro \circ f$ and
conversely $\ro = T\ovl f \circ \ro' \circ f^{-1}$. If $\ovl f$ is
an injective immersion, then we say that $(\nu',\ro')$ is an
anchored subbundle of $(\nu,\ro)$. Note that $\ro'$ is completely
determined by $\ro'=T (\ovl f)^{-1} \circ \ro \circ f$, which is
well defined since $\ovl f$ is an immersion. Assume that both
anchored bundles are linear. Then, we say that $f$ is a {\em
linear homomorphism} if $f: N' \to N$ is a linear bundle map. The
following commutative diagram represents an anchored bundle
morphism:

\begin{center}
\setlength{\unitlength}{1cm}
\begin{picture}(8,3.5)(.5,.5) \thicklines
\put(0.721,2.31){\vector(2,-3){0.968}}\put(0.360,2.37){\rmfamily
$N'$} \put(0.911,2.52){\vector(1,0){2.25}}
\put(3.35,2.25){\vector(-3,-4){1.05}} \put(1.73,0.466){\rmfamily
$M'$} \put(3.33,2.35){\rmfamily $TM'$} \put(0.7,1.42){\rmfamily
$\nu'$} \put(3.03,1.42){\rmfamily $\t_{M'}$}
\put(1.89,2.63){\rmfamily $\ro'$}
\put(2.50,0.678){\vector(3,1){3.5}} \put(4.26,1.46){\rmfamily
$\ovl f$} \put(0.890,2.65){\vector(3,1){3.58}}
\put(5.13,3.88){\vector(1,0){2.25}} \put(6.00,4.03){\rmfamily
$\ro$} \put(7.63,3.60){\vector(-3,-4){1.05}}
\put(6.05,1.93){\rmfamily $M$} \put(7.2,2.71){\rmfamily $\t_M$}
\put(5.02,3.69){\vector(2,-3){0.968}} \put(5.2,2.71){\rmfamily
$\nu$} \put(2.80,3.46){\rmfamily $f$} \put(5.7,3.4){\rmfamily
$T\ovl f$} \put(4.22,2.67){\vector(3,1){3.19}}
\put(7.57,3.69){\rmfamily $TM$} \put(4.66,3.73){\rmfamily $N$}
\end{picture}
\end{center}

\subsection{The foliation on anchored bundles} In this section we
need some elements of the theory of integrability of
distributions, developed by H.J. Sussmann \cite{Suss} (see also
\cite{Lib}). We first briefly recall the basic definitions and
main results on distributions, before applying them to anchored
vector bundles. We also use this section to fix some notations
regarding composite flows and concatenations of integral curves of
vector fields.

Consider a manifold $M$ and assume that $F$ is a differentiable
distribution on $M$, i.e. $F$ is a subset of $TM$ such that, for
any point $x \in M$, the fibre $F_x = F\cap T_xM$ is a linear
subspace of $T_xM$ and such that $F_x$ is spanned by a finite
number of vector fields in $F$ evaluated at $x$ (we say that $X
\in \vectorfields{M}$ is a vector field in $F$ if $X(y) \in F_y$,
for arbitrary $y\in M$).

The {\em rank} of the distribution $F$ at a point $x \in M$ is the
dimension of $F_x$. Note that, in the above definition, a
distribution need not have, in general, constant rank. If $F$ has
constant rank, we say that $F$ is a {\em regular distribution}.

A distribution is said to be completely integrable if, given any
$x \in M$, then there exists an immersed connected submanifold $i:
L \hookrightarrow M$ containing $x$ and such that $T_yL = F_y$,
for each $y \in L$. A submanifold $L$ satisfying the above
conditions, is called a leaf if it is maximal, in the sense that,
given any other submanifold $L'$, verifying the above conditions,
and which contains $L$ then $L' = L$. It can be proven that these
leaves are unique and determine a partition on $M$ which is called
the foliation induced by the completely integrable distribution.
Note that, by definition, the distribution $F$ has constant rank
on the points of the leaf $L$.

Assume that $\cal{F}$ is a family of vector fields on $M$, each
defined on an open subset of $M$. We say that $\cal{F}$ is
everywhere defined if, given any $x \in M$, there exists an
element $X$ of ${\cal F}$ containing $x$ in its domain. An
everywhere defined family of vector fields $\cal{F}$ generates a
distribution $F$ in the following way:
\[
F_x = \spn \{ X(x) \ | \ X \in \cal{F}, x \in \dom X\}.
\]
It is readily seen that $F$ is a differentiable distribution. H.J.
Sussmann has shown that, given of an everywhere defined family of
vector fields $\cal F$, generating a distribution $F$, one can
always construct the smallest completely integrable distribution
containing $F$. In order to discuss this construction, we need the
notion of a composite flow.

Assume that we have fixed an ordered $\ell$-tuple
$\cal{X}=(X_\ell,\ldots,X_1)$ of vector fields on $M$, and let us
represent the flow of $X_i$ by $\{\phi^i_t\}$.

The {\em composite flow} of $\cal{X}$ is the map
\[
\Phi: V \subset \R^\ell \times B \to B: ((t_\ell ,\ldots, t_1),x)
\mapsto \phi^\ell_{t_\ell} \circ \ldots \circ \phi^1_{t_1}(x),
\]
defined on some open subset $V$ of $\R^\ell\times B$. For brevity
we shall write $\Phi_T(x)$ in stead of $\Phi (T,x)$, where
$T=(t_\ell,\ldots,t_1)$. We shall sometimes refer to $T$ as the
{\em composite flow parameter}. For each fixed $T$, $\Phi_T$
determines a diffeomorphism from an open subset of $M$ (which may
be empty) to another open subset of $M$. It can be proven that, if
we fix a point $x \in M$, then the map $T'\mapsto \Phi_{T'}(x)$ is
smooth and defined on an open neighbourhood of $T$. For further
details on the domain of composite flows, we refer the reader to
\cite{Lib}.

Assume that we have fixed two composite flows: $\Phi$ of $\cal
X=(X_\ell,\ldots,X_1)$ and $\Psi$ of $\cal
Y=(Y_{\ell'},\ldots,Y_1)$. The {\em composition of} $\Phi$ and
$\Psi$ is the composite flow $\Psi \star \Phi$ of the
$\ell'+\ell$-tuple $(Y_{\ell'},\ldots,Y_1,X_\ell,\ldots,X_1)$.
Using these notations, it is easily seen that, for instance,
$\Phi$ equals $\phi^\ell \star \ldots \star \phi^1$. If $T$ is a
composite flow parameter for $\Phi$ and $T'$ for $\Psi$, then
$T'\star T=(T',T) \in \R^{\ell' + \ell}$ is a composite flow
parameter for $\Psi \star \Phi$.

The composite flow $\Phi$ of $\cal{X}= (X_\ell,\ldots,X_1)$ is
said to be {\em generated} by an everywhere defined family of
vector fields $\cal{F}$ if $\cal{X}$ is an ordered $\ell$-tuple of
elements of $\cal F$. Using all composite flows generated by $\cal
F$, we can define an equivalence relation on the points of $M$,
denoted by $\rel{}{}$.

\begin{defn}
Assume that $x, y \in M$. Then $\rel{x}{y}$ if there exists a
composite flow $\Phi$ generated by $\cal F$ and a composite flow
parameter $T$ such that $\Phi_T(x)=y$.
\end{defn}

It is easily seen that the relation $\rel{}{}$ is transitive (see
the above definition of the composition of composite flows) and
reflexive (take $T=(0,\ldots,0)$). If $\Phi$ is a composite flow
of $\cal{X}=(X_\ell,\ldots,X_1)$ and $\Phi_T(x)=y$ for some
$T=(t_\ell,\ldots,t_1)$, then the composite flow $\tilde \Phi$ of
$\tilde{\cal{X}}=(X_1,\ldots,X_\ell)$ and the composite flow
parameter $\tilde T=(-t_1,\ldots,-t_\ell)$ satisfy $\tilde
\Phi_{\tilde T}(y)=x$. Since $\tilde \Phi$ is generated by $\cal
F$, this makes the relation symmetric. Assume in the following
that the distribution $F$ is generated by the family $\cal F$.

\begin{thm}\label{thmsus}
The smallest completely integrable distribution $\widetilde{F}$
containing $F$ is the distribution generated by the everywhere
defined family $\widetilde{\cal{F}}$ containing all elements of
the form $\Phi^*_T Y$ where $Y \in \cal{F}$ and $\Phi$ is a
composite flow generated by $\cal F$.

The leaves of the distribution $\widetilde{F}$ are the equivalence
classes of the equivalence relation $\rel{}{}$.
\end{thm}

Consider the distribution $\widetilde F$ and let $[X,Y]$ denote
the Lie bracket of two vector fields in $\cal F$. It is easily
seen that $[X,Y]$ is a vector field in $\widetilde{F}$. Indeed,
let $\{\phi_t\}$ be the flow of $X$ and observe that $\phi_t^*Y$
is in $\widetilde{\cal F}$. Then, for each $x \in M$, the curve $t
\mapsto \phi^*_tY(x)$ is entirely contained in the linear space
$\widetilde F_x$, and so is its tangent vector:
\[
\left.\frac{d}{dt}\right|_0 \left(\phi^*_t Y(x)\right) =[X,Y](x).
\]
This reasoning can easily be extended to any finite number of
iterated Lie brackets of vector fields in $\cal F$. In fact, this
observation is rather important since it leads to an alternative
proof for the Theorem of Chow (see \cite{Suss}).

Assume that $\cal{X}=(X_\ell,\ldots,X_1)$ is an arbitrary finite
ordered family of vector fields, with composite flow $\Phi$. Fix a
value $(t_\ell,\ldots,t_1)$ of the composite flow parameter $T$.
The {\em concatenation of integral curves} of $\cal X$ trough $x
\in M$ is the piecewise smooth curve $\g: [a,a+|t_1|+\ldots
+|t_\ell|]\to M$ defined as follows, where $a_i = a+
\sum_{j=1}^i|t_j|$, ${\sign(t_i)} = \frac{t_i}{|t_i|}$ for $t_i
\neq 0$ and $\sign(0)=0$,
\[
\g(t) = \left\{\begin{array}{lll} \phi^1_{\sign(t_1)(t-a)}(x) &
\mbox{ for } & t \in [a,a_1]\\
\phi^2_{\sign(t_2)(t-a_1)}(\phi^1_{t_1}(x))&\mbox{ for } &
t\in\ ]a_1,a_2]\\
\ldots& & \\
\phi^\ell_{\sign(t_\ell)(t-a_{\ell-1})}(\ldots (\phi^1_{t_1}(x))
\ldots )&\mbox{ for } & t\in \
]a_{\ell-1},a_\ell],\end{array}\right.
\]
Note that, if $t \in ]a_{i-1},a_i[$ then $\dot \g(t) = \sign(t_i)
X_i(\g(t))$ and, hence, the restriction of $\g$ to $]a_{i-1},a_i[$
is an integral curve of $X_i$ if $t_i > 0$ (or $-X_i$ if $t_i
<0$). Note that $\g(a_\ell) = \Phi_T(x)$, i.e. the endpoint of
$\g$ coincides with the image of $x$ under the composite flow
$\Phi_T$. It is easily seen that in the specific case where $\cal
X$ is generated by $\cal F$, the concatenation of integral curves
of $\cal X$ through $x \in M$ is entirely contained in the leaf
$L_x$ through $x$.

Let us now proceed towards the construction of an everywhere
defined family of vector fields on $M$, given an anchored bundle
$(\nu,\ro)$ on $M$. Consider an arbitrary (local) section $\s$ of
$\nu$, i.e. $\s: M \to N$ is a smooth map with $(\nu \circ
\s)(x)=x$. Using the anchor map we can define the following vector
field on $M$: $\ro \circ \s$. Let $\cal D$ denote the set of all
vector fields of the form $\ro \circ \s$. Clearly, $\cal D$ is
everywhere defined and using the notations as described above, the
manifold $M$ is equipped with a distribution $D$ generated by
$\cal D$ (with $D = \im \ro$ if $(\nu,\ro)$ is linear) and the
smallest completely integrable distribution $\widetilde D$
containing $D$. The leaf on $M$ through $x$, induced by
$\widetilde D$, is denoted by $L_x$.

Consider the immersion $i:L_x \hookrightarrow M$, and let $\nu':N'
=L_x \times_{M} N\to L_x$ denote the pull-back bundle of $\nu$
under $i$, i.e. $(y,s) \in N'$ if $i(y) = \nu(s)$. Since $i$ is an
immersion, we can define an anchor map $\ro': N'\to TL_x$ as
follows: $T_yi(\ro'(y,s)) = \ro(s)$, given any $(y,s) \in N'$. The
projection $\pi_2:N' \to N$ of $N'$ onto the second factor,
determines an anchored bundle morphism, fibred over the immersion
$i$, i.e. $(\nu',\ro')$ is an anchored subbundle of $(\nu,\ro)$.
We shall call $(\nu',\ro')$ {\em the pull-back anchor bundle under
$i$}.

Before passing to the next section, we first give two examples of
an anchored bundle and the distribution induced by it. The first
example is taken from \cite{Mont}, where it was used in the
context of sub-Riemannian geometry to construct length-minimising
strictly abnormal extremals. The other example is taken from
\cite{Lib} and provides a non-trivial completely integrable
distribution on $\R^2$.

\begin{exmp} \label{voorbeeldmont}\textnormal{
Assume that $M=\R^3$ (we use cylindrical coordinates on $\R^3$),
and that $\nu: N = \R^3 \times \R^2 \to \R^3$ is a trivial bundle
over $M$. Consider the following two vector fields on $M$: $X_1 =
\fpd{}{r}$ and $X_2 = \fpd{}{\theta} - p(r) \fpd{}{z}$, where
$p(r)$ is a function on $\R$ with a single non degenerate maximum
at $r=1$, i.e. $p$ satisfies:
\[
\left. \frac{d}{dr} p(r) \right|_{r=1} = 0 \quad \hbox{\rm and}
\quad \left. \frac{d^2}{dr^2} p(r) \right|_{r=1} < 0.
\]
Such a function can always be constructed (take, for instance,
$p(r) = \half r^2 - \frac{1}{4}r^4$). Let $\ro$ denote the map
defined by $\ro(x,u^1,u^2) = u^1X^1(x) + u^2X^2(x)$, with
$x=(r,\theta,z) \in M$. It is easily seen that $(\nu,\ro)$ is a
linear anchored bundle. The flows of $X_1,X_2$ are denoted by
$\{\phi_t\}$, $\{\psi_t\}$, respectively. In particular, we have
$\phi_t(r,\theta,z)= (t+r,\theta,z)$, $\psi_t(r,\theta,z)=
(r,\theta +t,z-p(r)t)$. The foliation induced by $\im \ro$ is
trivial. Indeed, all iterated Lie brackets of the two vector
fields $X^1$ and $X^2$ span the total tangent space at each point,
implying that $\widetilde D = TM$ and $M$ itself is the only
leaf.}\end{exmp}

\begin{exmp} \label{voorbeeld2}\textnormal{
Let $M= \R^2$ and let $N=M\times \R^2$, with $\ro(x,u^1,u^2) =
u^1X(x) + u^2 Y(x)$, where \[ X=\fpd{}{x} \mbox{ and }
Y=y\fpd{}{y}.\] The distribution $F$ on $M$ defined by $F=\im \ro
$ satisfies $F= \widetilde F$, since $[X,Y]=0$, i.e. $F$ is
completely integrable. The two 2-dimensional submanifolds
$\{y<0\}$, $\{y>0\}$ and the 1-dimensional submanifold $\{y=0\}$
are the leaves of the foliation on $M$. We use this example to
show that Lemma \ref{lem0} in the following section is
non-trivial. We shall construct a curve, which is tangent to $F$,
i.e. has tangent vector everywhere contained in $F$, but, the
curve itself is not entirely contained in a single leaf. Indeed,
consider $\tilde c: \R \to M : t \mapsto (t,t^3)$. It is readily
seen that $\dot{\tilde{c}}(t) = X(\tilde{c}(t)) + 3 t^{-1}
Y(\tilde c(t)) \in F_{\tilde c(t)}$ for $t\neq 0$ and $\dot{\tilde
c}(0) = X(0,0) \in F_x$. However $\tilde c$ passes through the
three leaves of $F$.}
\end{exmp}

\subsection{$\ro$-admissible curves}\label{skrommen}

We introduce here the notion of a {\em $\ro$-admissible curve}. By
a smooth curve in a manifold $M$ we will always mean a
$C^{\scriptscriptstyle \infty}$ map $c: I \rightarrow M$, where $I
\subseteq \R$ may be either an open or a closed (compact)
interval. In the latter case, the denominations ``path" or ``arc"
are also frequently used in the literature but, for simplicity, we
will make no distinction in terminology between both cases. For a
curve defined on a closed interval, say $[a,b]$, it is tacitly
assumed that it admits a smooth extension to an open interval
containing $[a,b]$. Fix an anchored bundle $(\nu,\ro)$ on $M$.

\begin{defn}
Let $c: [a,b] \to N$ denote a smooth curve in $N$, and let $\tilde
c = \nu \circ c$ denote the projected curve in $M$, called the
base curve of $c$. Then $c$ is called a {\em smooth
$\ro$-admissible curve} if $\ro \circ c = \dot{\tilde c}$.
\end{defn}

Local coordinates on $M$ will be denoted by $(q^i)$ and
corresponding bundle adapted coordinates on $N$ by $(q^i,u^{a})$,
with $i=1,\ldots,n$ and $a=1,\ldots,k$, where $k$ is the dimension
of the typical fibre of $N$. If we write the bundle map $\ro$
locally as
\begin{equation} \label{rolocal}
\ro(q^i,u^a) = \g^i(q^j,u^a)\fpd{}{q^i}
\end{equation}
Then a smooth $\ro$-admissible curve $c(t) = (q^i(t),u^a(t))$
locally satisfies
\[
\g^i(q^j(t),u^a(t)) = \dot q^i(t).
\]
In order to introduce a suitable concept of ``leafwise holonomy''
in the framework of principal $\rho$-lifts, it turns out that the
class of $\ro$-admissible curves in $N$ should be further extended
to curves admitting (a finite number of) discontinuities in the
form of certain `jumps' in the fibres of $N$, such that the
corresponding base curve is piecewise smooth. In order to define
these ``piecewise'' $\ro$-admissible curves we first consider the
composition of smooth $\ro$-admissible curves.

The {\em composition} of a finite number of, say $\ell$, smooth
$\ro$-admissible curves $c_i:[a_{i-1},a_i] \to N$ for
$i=1,\ldots,\ell$, satisfying the conditions $\tilde c_i(a_{i}) =
\tilde c_{i+1}(a_i)$ for $i=1,\ldots,\ell-1$, is the map $c_\ell
\cdot \ldots \cdot c_1:[a_0,a_\ell] \to N$ defined by
\begin{equation}\label{compositie}
(c_\ell\cdot \ldots \cdot c_1) (t) = \left\{ \begin{array}{ll}
c_1(t)& \quad t \in [a_0,a_1],\\
\ldots\\ c_\ell(t) & \quad t \in \
]a_{\ell-1},a_\ell].\end{array}\right.
\end{equation}
Note that the base curve of $c_\ell\cdot\ldots \cdot c_1$ is a
piecewise smooth curve. However, in general
$c_\ell\cdot\ldots\cdot c_1$ is discontinuous at $t=a_i$,
$i=1,\ldots,\ell-1$. The composition $c=c_\ell \cdot \ldots \cdot
c_1$ is called {\em a piecewise $\ro$-admissible curve}, or simply
a {\em $\ro$-admissible curve}. We now proceed towards the
following important result, saying that the base curve of a
$\ro$-admissible curve is always entirely contained in a leaf of
the foliation on $M$, induced by the everywhere defined family of
vector fields $\cal D$ on $M$ (see the previous section).

\begin{lem}\label{lem0} The base curve $\tilde c$ of a $\ro$-admissible
curve $c:[a,b]\to N$ is entirely contained in the leaf $L_x$, with
$x=\tilde c(a)$.
\end{lem}
\begin{pf}
It is sufficient to prove this result for $c$ smooth. For any
point $y\in M$, consider a coordinate neighbourhood $U$ of $y$
with coordinates $(q^1,\ldots,q^n)$, adapted to the foliation
induced by $\cal D$, such that: $(1)$ if $q^{p+1}(z),\ldots,q^n(z)
=0$ then $z \in L_y$ and $(2)$ the coordinate functions
$q^1,\ldots,q^p$ determine local coordinates on the leaf $L_y$
(this is always possible since $L_y$ is an immersed submanifold).
Upon restricting $U$ to a smaller subset, if necessary, we may
always assume, in addition, that the fibre bundle $\nu$ is trivial
over $U$, and we denote the adapted bundle coordinates by
$(q^i,u^a)$, for $i=1,\ldots,n$ and $a=1,\ldots,k$. In the
following we only consider such coordinate charts. Recall the
definition of the pull-back anchored bundle $(\nu',\ro')$ under
$i:L_y \hookrightarrow M$. Note that
$(q^1,\ldots,q^p,u^1,\ldots,u^k)$ is a bundle adapted coordinate
chart on $N'$.

Fix a coordinate chart (in the sense specified above) containing
the point $x= \tilde c(a)$ and assume that $c$ is written in these
coordinates as $(\tilde c^i(t),u^a(t))$. Let $\tilde d$ denote the
solution in $L_y$ of the following differential equation, in the
anchored bundle $(\nu',\ro')$:
\[
\dot{\tilde d}(t)=\ro'(\tilde d^1(t),\ldots,\tilde
d^p(t),u^1(t),\ldots,u^k(t)),
\]
with initial condition $\tilde d(a)=x$. From standard arguments we
know that $\tilde d$ is defined on some interval, say $[a,t+\e[$,
with $\e>0$.

Consider the curve $\tilde d'=i\circ \tilde d:[a,t+\e[\to M$ in
$M$, through $y$ at time $t$. Then we have, by uniqueness of
solutions to differential equations, that $\tilde d' = \tilde
c|_{[a,t+\e[}$, since the curves $\tilde d'$ and $\tilde c$ both
solve $\dot{\tilde c} = \ro(\tilde c^j,u^a)$. Indeed, for $\tilde
c$ this is trivially satisfied and for $\tilde d'$ we have
\[
\dot{\tilde d'}(t) = Ti (\ro'(\tilde d(t),c(t))) = \ro(\tilde
d^i,0,u^a(t')).
\]
Therefore we conclude that $\tilde c|_{[a,t+\e[}$ is contained in
the leaf $L_x$, since by making use of the coordinate system, we
have $\tilde c^i(t) = 0$ for $t \in [a,t+\e[$ and
$i=p+1,\ldots,n$. Taking the limit from the left at $t=a+\e$, we
obtain that $\tilde c^i(a+\e)=0$ for $i=p+1,\ldots,n$, or $\tilde
c(a+\e) \in L_x$. We can repeat the above reasoning for the curve
$c|_{[a+\e,b]}$, i.e. we start from the point $\tilde c(a+\e)$ in
stead of the point $x$. We thus obtain that $\tilde c(t) \in L_x$
for all $t\in[a,a+\e+\e']$ for some $\e' >0$. Continuing this way,
we eventually obtain that the entire curve $\tilde c$ is contained
in $L_x$, concluding the proof. \qed
\end{pf}

It can be seen that the curve $\tilde c$ constructed in Example
\ref{voorbeeld2} does not contradict the previous lemma although
$\tilde c$ is a curve tangent to the distribution $\im \ro$.
Indeed, $\tilde c$ can not be written as the base curve of a
$\ro$-admissible curve, since, at $t=0$ a singularity is
encountered.

Consider two anchored bundles $(\nu',\ro')$ and $(\nu,\ro)$, and a
anchored bundle morphism $f$ between them, i.e. $f:N' \to N$
fibred over $\ovl f:M'\to M$. Let $c'$ denote a $\ro'$-admissible
curve. Consider the curve $c=f\circ c'$ in $N$, and let $\tilde
c$, resp. $\tilde c'$, denote the base curve of $c$, resp. $c'$.
Then, we have that $c$ is $\ro$-admissible, since
\[
\ro \circ c = \ro \circ f \circ c' = T \ovl f \circ \ro' \circ c'
= T \ovl f \circ \dot{\tilde c}' = \dot{\tilde c}.
\]

Let $c$ denote a $\ro$-admissible curve. If $x=\tilde c(a)$ and
$y=\tilde c(b)$, then we say that {\em $c$ takes $x$ to $y$}, and
we write $x \stackrel{c}{\to} y$ (shortly $x\to y$ if we do not
want to mention the $\ro$-admissible curve explicitly). The
relation $\to$ on $M$ is transitive, and is preserved by an
anchored bundle morphism, i.e. if $x'\to y'$ then $\ovl f(x')\to
\ovl f(y')$ for $x',y' \in M'$. The set of points $y$ such that $x
\to y$ for some fixed $x$ is denoted by $R_x$ and is called {\em
the set of reachable points from $x$}. Until now, we have proven
that the base curve of a $\ro$-admissible curves is contained in a
leaf $L_x$ of the foliation on $M$, i.e. $R_x \subset L_x$. It is
interesting to wonder if every point in $L_x$ can be reached from
$x$ following a $\ro$-admissible curve. In general this is not the
case. However, if we consider the composition of $\ro$- and
$(-\ro)$-admissible curves, then every point in $L_x$ can be
reached.

\begin{defn}
Given an anchored bundle $(\nu,\ro)$. The {\em inverse anchored
bundle} is defined as $(\nu,-\ro)$, where $-\ro: N \to TM:
s\mapsto -\ro(s)$.
\end{defn}

An anchored bundle $(\nu,\ro)$ is related to its inverse in the
following way. Assume that $c$ is a $\ro$-admissible curve taking
$x$ to $y$, i.e. $x \stackrel{c}{\to} y$. Then the curve $c^*:
[a,b]\to N: t \mapsto c((b-t)+a)$ is $(-\ro)$-admissible and takes
$y$ to $x$. We shall call this curve {\em the $(-\ro)$-admissible
curve associated with $c$}, or simply the {\em reverse of $c$}.
Note that, using these notations, $(c^*)^*=c$. If we write, the
relation on $M$ induced by the inverse anchored bundle as $\ot$,
we have the following equivalence:
\[ x \stackrel{c}{\to} y \mbox{ iff } y \stackrel{c^*}{\ot} x.\]

Note that the family of vector fields on $M$ defined by the
inverse anchored bundle equals $-\cal D=\{-\ro\circ\s \ | \ \s \in
\G(\nu)\}$, and therefore produce the same distribution $D$ and
the same foliation as $\cal D$. We now consider the composition of
$\ro$- and $(-\ro)$-admissible curves. Thus, assume that we have
$\ell$ curves $c_i: [a_{i-1},a_i]\to N$ for $i=1,\ldots,\ell$ such
that $c_{i-1}(a_{i-1})=c_i(a_{i-1})$ and such that $c_i$ is either
$\ro$-admissible or $(-\ro)$-admissible. The composition of the
curves $c_i$ (defined as in Equation \ref{compositie})
$c=c_\ell\cdot\ldots\cdot c_1$ is called a {\em
$\pm\ro$-admissible curve}.

The projection $\tilde c$ of $c$ onto $M$ is a piecewise smooth
curve which is called {\em the base curve of the
$\pm\ro$-admissible curve}. If $\tilde c(a_0)=x$ and $\tilde
c(a_\ell)=y$ we say that the $\pm\ro$-admissible curve takes $x$
to $y$. Note that, in this case, the $\pm\ro$-admissible curve
$c^*$, defined by $c^*= (c_{1})^*\cdot\ldots\cdot (c_{\ell})^*$,
takes $y$ to $x$.

We thus obtain an alternative characterisation of the leaves of
the foliation generated by the anchored bundle $(\nu,\ro)$.

\begin{thm} \label{thmkarleafs}
We have that $\rel{x}{y}$, or $y \in L_x$, iff there exists a
$\pm\ro$-admissible curve taking $x$ to $y$.
\end{thm}
\begin{pf}
The `if'-part of the proof follows straightforwardly from Lemma
\ref{lem0}. The `only if'-part is proven by the following
reasoning. Assume that $y\in L_x$ and consider a composite flow
$\Phi$ of $\cal{X}=(X_\ell,\ldots,X_1)$, with $X_i= \ro \circ
\s_i$ and $\s_i \in \G(\nu)$ ($\Phi$ is generated by $\cal{D}$)
such that $\Phi_T(x) =y$. Consider the following curves,
\[
c_i:[a_{i-1},a_i] \to N: t \mapsto \s_i \circ \g|_{[a_{i-1},a_i]},
\]
where $\g$ is the concatenation associated with $\cal X$ and $T$
through $x$ (where we have used the notations from the preceding
section). It is easily seen that $c_i$ is $\ro$-admissible if
$\sign (t_i) > 0$, and $(-\ro)$-admissible if $\sign (t_i) < 0$.
If we put $c = c_{\ell}\cdot \ldots\cdot c_1$, then $c$ takes $x$
to $y$ and is $\pm \ro$-admissible. \qed
\end{pf}

The proof of the following theorem now easily follow from Theorem
\ref{thmkarleafs}. Note that any anchored bundle morphism $f$
between $(\nu',\ro')$ and $(\nu,\ro)$, which is fibred over $\ovl
f:M'\to M$, is also a morphism of the corresponding inverted
anchored bundles, i.e. $f:(\nu',-\ro')\to (\nu,-\ro)$. This
implies that, if $x' \ot y'$ then $\ovl f(x') \ot \ovl f(y')$, for
$x',y' \in M'$.

\begin{thm}
Let $f$ denote a morphism between $(\nu',\ro')$ and $(\nu,\ro)$,
fibred over $\ovl f:M'\to M$. Then $\ovl f(L_{x'}) \subset L_{\ovl
f(x')}$. If $(\nu',\ro')$ is the pull-back bundle along $i:L_x
\hookrightarrow M$ and $f=\pi_2$, then $i(L_{x}) = L_{i(x)}$.
\end{thm}

It is interesting to consider the special case of {\em linear}
anchored bundles.

\begin{thm}\label{thmkrommenenleaf}
Let $(\nu,\ro)$ denote a linear anchored bundle on $M$ and take
any $x, y \in M$. Then $y \in L_x$ or $\rel{x}{y}$ iff there
exists a $\ro$-admissible curve that takes $x$ to $y$, i.e. we
have $R_x=L_x$.
\end{thm}

This theorem follows from the fact that, given a linear anchored
bundle, then $x \to y$ iff $y\to x$. Indeed, assume that
$c:[a,b]\to N$ is a $\ro$-admissible curve taking $x$ to $y$. Then
the curve $c^{-1}:[a,b]\to N: t \mapsto -c((b-t)+a)$ is also
$\ro$-admissible and takes $y$ to $x$. Note that $c^{-1}=-c^*$.
The curve $c^{-1}$ is called {\em the inverse of c}. In
particular, the base curve of a $\pm\ro$-admissible curve is the
base curve of a $\ro$-admissible curve on a linear anchored
bundle, which proves the above theorem. Let $c:[a,b]\to N$ denote
a smooth $\ro$-admissible curve. We now prove that any
``reparameterisation'' of $\tilde c$ is again a base curve of a
$\ro$-admissible curve. Assume that $\phi: [a,b] \to [c,d]$ is a
diffeomorphism satisfying $\phi(a)=c$ and $\phi(b)=d$. Consider
the following curve $c':[c,d]\to N$ defined by
\[
c'(s) = \frac{d\phi^{-1}}{ds}(s) c(\phi^{-1}(s)).
\]
From elementary calculations, it is easily seen that $c'$ is
$\ro$-admissible, and that the base curve equals $\tilde
c(\phi^{-1}(s))$, i.e. the reparameterisation of $c$. Note that
the above definitions are only valid if $(\nu,\ro)$ is a linear
anchored bundle.

\subsection{$\pm\ro$-admissible loops}

Consider a point $x \in M$ and let $C(x,N)$ denote the set of all
$\pm\ro$-admissible curves taking $x$ to itself. Elements of
$C(x,N)$ are called, with some abuse of terminology, {\em
$\pm\ro$-admissible loops with base point $x$}. Indeed, in general
a $\pm\ro$-admissible loop need not be continuous, nor closed.

Let $\pi_1(x,M)$ denote the first homotopy group of $M$ with
reference point $x$ and consider the map $C(x,N) \to \pi_1(x,M)$,
associating to the base curve of a $\pm\ro$-admissible loop $c$,
its homotopy class in $\pi_1(x,M)$, i.e. if $\tilde c$ is the base
curve of $c=c_\ell\cdot\ldots\cdot c_1 \in C(x,N)$, then $\tilde
c$ is mapped onto $[\tilde c]$. It is easily seen that the image
of $C(x,N)$ determines a subgroup of $\pi_1(x,M)$, which is
denoted by $\pi_1^N(x,M)$. Indeed, assume that
$c=c_\ell\cdot\ldots\cdot c_{1}$ and $d=d_{\ell'}\cdot\ldots\cdot
d_{1}$ are elements of $C(x,N)$, with homotopy classes $[\tilde
c]$ and $[\tilde d]$ in $\pi_1(x,M)$. Then, the product $[\tilde
c] \cdot [\tilde d]$ in $\pi_1(x,M)$ is the homotopy class of the
base curve of
\[
c_\ell\cdot\ldots\cdot c_1\cdot d_{\ell'}\cdot \ldots\cdot d_1.
\]
On the other hand, if $c=c_\ell\cdot\ldots\cdot c_1$ is a
$\pm\ro$-admissible loop with base point $x$, then the curve
$c^*=(c_1)^*\cdot\ldots\cdot (c_\ell)^*$ is also contained in
$C(x,N)$, and the homotopy class of the base curve of $c^*$ is
precisely the inverse $[\tilde c]^{-1}$ of $[\tilde c]$.
Therefore, the $\pm\ro$-admissible loops generate a subgroup of
$\pi_1(x,M)$ which is denoted by $\pi_1^N(x,M)$. Note that, if
$(\nu,\ro)$ is linear, then $\pi_1^N(x,M)$ is generated by the set
of $\ro$-admissible loops with base point $x$, i.e.
$\ro$-admissible curves taking $x$ to itself.

We now elaborate on how the above defined structures on anchored
bundles behave under homomorphisms. From Section \ref{skrommen},
we already now that $\pm\ro$-admissible curves are preserved under
anchored bundle morphisms. Similarly, $\pm\ro$-admissible loops
are preserved, taking us to a group morphism between the
corresponding subgroups of the first fundamental group of the base
manifolds. More precisely, assume that $f$ denotes a homomorphism
between two anchored bundles $(\nu',\ro')$ and $(\nu,\ro)$, fibred
over $\ovl f$. Then, if $[\ovl f]$ denotes the corresponding group
morphism from $\pi_1(x',M')$ to $\pi_1(\ovl f(x'),M)$, we have
that $[\ovl f]$ can be restricted to a morphism from
$\pi_1^{N'}(x',M')$ to $\pi_1^N(\ovl f(x'),M)$.

Consider the pull-back setting under $i:L_x \hookrightarrow M$,
and let $\pi_2: N'=i^*N \to N$ denote the associated anchored
bundle morphism. From the above, we now that $[i]$ maps the
subgroup $\pi_1^{N'}(y,L_x)$ of $\pi_1(L_x)$ to the subgroup
$\pi_1^N(y,M)$ $\pi_1(y,M)$ (note that $L_x$ is connected,
allowing us to omit the reference point in the first homotopy
group of $L_x$). We now prove that $[i]:\pi_1^{N'}(y,L_x) \to
\pi_1^N(y,M)$ is onto. Consider an arbitrary element of
$\pi_1^N(y,M)$ associated with some $c \in C(y,N)$. From Lemma
\ref{lem0} we know that $\tilde c$ is contained in the leaf
$L_y=L_x$, which in turn implies that there exists a
$\pm\ro$-admissible loop $c' \in C(y,N')$ such that $\pi_2\circ
\tilde c'=c$. In particular, we have that $[i]([\tilde c']) =
[\tilde c]$, and, hence, $[i]$ is onto, when restricted to
$\pi_1^{N'}(y,L_x)$.

Moreover, from Theorem \ref{thmkrommenenleaf}, we know that any
two points in $L_x$ can be connected by the base curve of a
$\pm\ro$-admissible curve. This implies, using standard arguments,
that we can omit the reference point in $\pi_1^{N'}(y,L_x)$, and
from now on, we use the notation $\pi_1^{N}(L_x)$ for
$\pi_1^{N'}(y,L_x)$. Similarly, we write $\pi_1^N(L_x,M)$ for
$\pi_1^N(x,M)$.

\section{Principal $\ro$-lifts}\label{sectieprincipalcon}

Let us first briefly recall the notion of a connection over a
bundle map in the context of principal fibre bundles and describe
some elementary properties. For further details we refer to
\cite{algcon}. Let $(\nu,\ro)$ denote an anchored bundle on $M$
and let $\pi: P\to M$ denote a principal fibre bundle with
structure group $G$. Consider the pull-back bundle $\tilde
\pi_1:\pi^*N\to P$ and let $\tilde \pi_2:\pi^*N\to N$ denote the
projection onto the second factor.

\begin{defn} \label{def1}A {\em principal lift over the bundle map} $\ro$, simply a
{\em principal $\ro$-lift}, is a bundle map $h : \pi^*N \to TP$,
fibred over the identity on $P$, such that in addition the
following conditions are satisfied for any $(u,s) \in \pi^*N$:
\begin{enumerate}
\item $TR_g(h(u,s))=h(ug,s)$, and
\item $T\pi\circ h = \ro\circ\tilde \pi_2$.
\end{enumerate}
If $(\nu,\ro)$ is a linear anchor bundle, the bundle $\tilde
\pi_1:\pi^*N\to P$ can be given linear structure. In this case, we
say that a principal $\ro$-lift $h$ is a {\em principal
$\ro$-connection} if $h:\pi^*N\to TP$ is, in addition, a linear
bundle morphism from $\tilde \pi_1$ to $\t_P$
\end{defn}

It is easily seen from the definition of a principal $\ro$-lift
$h$ that $(\tilde \pi_1,h)$ determines an anchored bundle and that
the projection $\tilde \pi_2: \pi^*N \to N$, which a bundle
morphism fibred over $\pi: P\to M$, determines an anchored bundle
morphism between $(\tilde \pi_1,h)$ and $(\nu,\ro)$. Moreover, if
$h$ is a $\ro$-connection, we have that $(\tilde\pi_1,h)$ is a
linear anchored bundle and that $\tilde \pi_2$ is a linear
anchored bundle morphism. The situation is illustrated by the
following diagram:

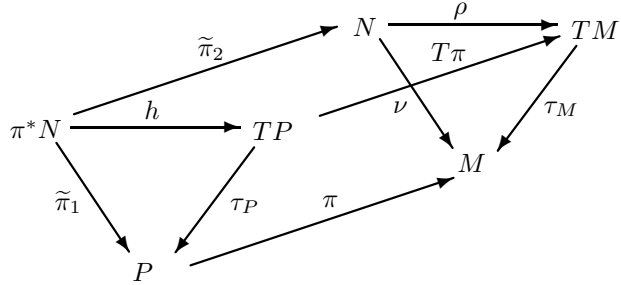
\begin{center}
\setlength{\unitlength}{1cm}
\begin{picture}(8,3.9)(.5,0) \thicklines
\put(0.721,2.31){\vector(2,-3){0.968}}\put(0.1,2.37){\rmfamily
$\pi^*N$} \put(0.911,2.52){\vector(1,0){2.25}}
\put(3.35,2.25){\vector(-3,-4){1.05}} \put(1.73,0.466){\rmfamily
$P$} \put(3.33,2.35){\rmfamily $TP$} \put(0.7,1.42){\rmfamily
$\tilde \pi_1$} \put(3.03,1.42){\rmfamily $\t_{P}$}
\put(1.89,2.63){\rmfamily $h$} \put(2.50,0.678){\vector(3,1){3.5}}
\put(4.26,1.46){\rmfamily $\pi$}
\put(0.95,2.69){\vector(3,1){3.5}}
\put(5.13,3.88){\vector(1,0){2.25}} \put(6.00,4.03){\rmfamily
$\ro$} \put(7.63,3.60){\vector(-3,-4){1.05}}
\put(6.05,1.93){\rmfamily $M$} \put(7.2,2.71){\rmfamily $\t_M$}
\put(5.02,3.69){\vector(2,-3){0.968}} \put(5.2,2.71){\rmfamily
$\nu$} \put(2.6,3.46){\rmfamily $\tilde\pi_2$}
\put(5.7,3.4){\rmfamily $T\pi$}
\put(4.22,2.67){\vector(3,1){3.19}} \put(7.57,3.69){\rmfamily
$TM$} \put(4.66,3.73){\rmfamily $N$}
\end{picture}
\end{center}

We will now apply the tools from the previous section to the study
of principal $\ro$-lifts. We first fix some notations and make
some preliminary comments. The everywhere defined family of vector
fields on $P$ generated by $(\tilde \pi_1,h)$ is denoted by $\cal
Q$, and correspondingly, the distribution on $P$ generated by
$\cal Q$ is denoted by $Q$. We refrain from calling $Q$ a
horizontal distribution since for arbitrary $u \in P$ it may be
that ${Q}_u$ has non-zero intersection with the distribution of
vertical tangent vectors $V_u\pi = \ker T\pi$. Moreover, in
general ${\cal Q}_u + V_u\pi \neq T_uE$, i.e.\ ${\cal Q}_u$ and
$V_u\pi$ do not necessarily span the full tangent space $T_uP$.
The smallest integrable distribution containing $Q$ is denoted by
denoted by $\tilde Q$. The leaf of $\widetilde Q$ through an
arbitrary point $u \in P$ is written as $H(u)$. The principal
$\ro$-lift $h$ can be used to lift several kinds of objects from
the anchored bundle $(\nu,\ro)$ on $M$ to the anchored bundle
$(\tilde\pi_1,h)$ on $P$. For instance, given any (local) section
$s$ of $\nu$, we can define a mapping $s^h: P \rightarrow TP$ by
\begin{equation}\label{lift}
s^h(u) = h(u, s(\pi(u))).
\end{equation}
It is seen that, by construction, $s^h$ is smooth and verifies
$\tau_P(s^h(u)) = u$, i.e.\ $s^h$ is a (local) vector field on
$P$, called the {\em lift of the section $s$ with respect to $h$},
or simply {\em the lift of $s$} if no confusion can arise. Let us
denote by $\cal D^h$ the everywhere defined family of vector
fields on $P$ defined by the lift of (local) sections of $\nu$.

Next, we recall some definitions and results on principal fibre
bundles and principal connections from \cite{koba} since they will
be used extensively in the following sections. Let $\pi: P\to M$
denote a principle fibre bundle with structure group $G$. The Lie
algebra of $G$ is denoted by $\alg{g}$.

Consider the smooth map $\s_u:G \to P$ for each $u \in P$ defined
by $\s_u(g)=ug$. Then we have $T_e\s_u: \alg{g} \to V_u\pi$ and
$TR_g \circ T_e\s_u = T_e\s_{ug} \circ Ad_{g^{-1}}$, where $Ad_h:
\alg{g}\to \alg{g}$ denotes the adjoint action of $G$ on its Lie
algebra. Given any $A \in \alg{g}$, let $\s(A)$ denote the
vertical vector field on $P$ defined by $\s(A)(u) = T_e\s_u(A)$.
It is easily seen that $(R_g)_*\s(A) = \s(Ad_{g^{-1}}A)$.

A standard principal connection on $P$ is defined by {\em a
connection form} $\w$ on $P$, i.e. $\w$ is a $\alg{g}$-valued one
form on $P$ satisfying the following two conditions: $(1)$ for any
$A \in \alg{g}$, $\w(\s(A)) = A$, and $(2)$ for any $g \in G$,
$R_g^* \w = Ad_{g^{-1}} \cdot \w$. It is well known that $\w$ is
equivalently defined by a {\em horizontal lift} $h^\w: P \times_M
TM \to TP$, where $h^\w$ and $\w$ are related in the following
way: $h^\w(u,X) = \tilde X - T_e\s_u(\w(\tilde X))$ for any
$\tilde X \in T_uP$ satisfying $T\pi (\tilde X)=X$. From $(2)$ it
follows that $h^\w$ is right invariant, i.e.
$TR_g(h^\w(u,X))=h^\w(ug,X)$ for $X\in T_{\pi(u)}M$ and $g\in G$
arbitrary. For the sake of completeness, we mention that,
equivalently, a principal connection can be defined by the right
invariant distribution spanned by the image of $h^\w$, determining
a direct decomposition of $TP$, i.e. if $\im h^\w = H\pi$, then
$TP = H\pi \oplus V\pi$.

Before starting our study of principal $\ro$-lifts, we state the
following lemma, taking from \cite[p 69]{koba}.

\begin{lem} \label{lemkromingroup}
Let $G$ be a Lie group and $\alg{g}$ its Lie algebra. Let $Y_t$,
for $a \le t \le b$, define a continuous curve in $\alg{g}$. Then
there exists in $G$ a unique curve $g_t$ of class
$C^{\scriptscriptstyle 1}$ such that $g(a)= e$ and $\dot g_t
g_t^{-1} =Y_t$ for $a\le t\le b$.
\end{lem}

Let us now return to the general treatment of principal
$\ro$-lifts. Let $(\nu,\ro)$ denote an anchored bundle on $M$ and
let $P$ denote a principal fibre bundle on $M$ with structure
group $G$.

Fix a standard principal connection $\w$ on $P$. In the following
we will use the connection form $\w$ in order to obtain an
alternative description for a principal $\ro$-lift $h$. This
alternative description will allow us to easily derive some
properties of lifts of $\ro$-admissible curves with respect to $h$
(see below) using the theory of standard connections. Thus, let
$h$ be a given principal $\ro$-lift and consider the map $\chi :
\pi^*N \to \alg{g}$ defined by $\chi(u,s)=\w(h(u,s))$ for any
$(u,s) \in \pi^*N$. Note that the following relation holds
$\chi(ug,s) = Ad_{g^{-1}}\cdot \chi(u,s)$ and that $h(u,s) =
T_e\s_u(\chi(u,s)) + h^\w(u,\ro(s))$. We shall sometimes refer to
$\chi$ as {\em the coefficient of $h$ with respect to $\w$}. The
pair $(\w,\chi)$ determines uniquely the principal $\ro$-lift $h$,
in the following way. Given any connection form $\w$ on $P$ and a
map $\chi:\pi^*N \to \alg{g}$, such that $\chi$ transforms under
the right action in the following way: $\chi(ug,s)=
Ad_{g^{-1}}\cdot \chi(u,s)$, then the map $h:\pi^*N \to TP$,
defined by $h(u,s) = h^\w(u,s)+T_e\s_u(\chi(u,s))$, determines a
principal $\ro$-lift. Note that the coefficient of $h$ with
respect to $\w$ is precisely $\chi$.

\begin{thm}
Given any principal $\ro$-lift $h$, then the following properties
hold:
\begin{enumerate}
\item \label{prindistrigen}the family $\cal{D}^h$ generates the distribution $Q$,
and, hence the integrable distribution $\widetilde Q$;
\item \label{prinprojcurves}any $h$-admissible curve is mapped by $\tilde \pi_2$ onto a
$\ro$-admissible curve;
\item \label{prinliftcurves}given any $\ro$-admissible curve $c$ taking $x$ to $y$ and a point
$u\in P_{x}$, then there exists a unique $h$-admissible curve
projecting onto $c$ by $\tilde \pi_2$ and whose base curve in $P$
passes through $u$.
\end{enumerate}
\end{thm}
\begin{pf} Properties \ref{prindistrigen} and \ref{prinprojcurves} are
trivial. In order to prove \ref{prinliftcurves}, we fix a
principal connection $\w$ and consider the coefficient $\chi$ of
$h$ with respect to $\w$.

First we prove that, given any $\ro$-admissible curve $c:[a,b]\to
N$ with base curve $\tilde c$, then there always exists a
$h$-admissible curve whose base curve passes through $u \in
P_{\tilde c(a)}$ at $t=a$. First, consider the horizontal lift of
$\tilde c$ with respect to the principal connection $\w$, i.e.
$\tilde d^\w(t)$ is the unique curve satisfying $\dot d^\w(t) =
h^\w(d^\w(t), \dot{\tilde c}(t))$ and $d^\w(a)=u$. Let $g(t)$
denote the curve in $G$ satisfying the equation $TR_{g(t)^{-1}}
\dot g(t) = \chi(d^\w(t),c(t))$ and $g(a) = e$, with $e$ the unit
element of $G$. From \ref{lemkromingroup} the curve $g(t)$ always
exists and is unique. We now prove that $(d(t),c(t)) \in \pi^*N$,
with $d(t)=d^\w(t)g(t)$, is a $h$-admissible curve. Indeed, we
find that:
\begin{eqnarray*}
\dot d(t) &=& TR_{g(t)} (\dot d^\w(t)) + T_e\s_{d^\w(t)}(\dot g(t)),\\
&=& TR_{g(t)} h^\w(d^\w(t),\tilde c(t)) +
T_e\s_{d(t)}(TL_{g^{-1}(t)} \cdot \dot g(t)),\\
&=& h^\w(d(t), \tilde c(t)) + T_e\s_{d(t)} (Ad_{g^{-1}(t)} \cdot
\chi(d^\w(t),c(t))).
\end{eqnarray*}
From the definition of $\chi$, the tangent vector $\dot d(t)$
equals the desired vector $h(d(t),c(t))$. Clearly $(d(t),c(t))$
projects onto $c(t)$ and its base curve $d(t)$ passes through $u$
at $t=a$. It easily follows that $d(t)$ is uniquely determined by
these conditions, since it satisfies a first order differential
equation, i.e. $\dot d(t)=h(d(t),c(t))$, with given initial
condition $d(a)=u$. \qed
\end{pf}

{\bf $h$-Displacement and holonomy}

Using the notations from the above theorem, we have that $d(t)$ is
uniquely determined from the $\ro$-admissible curve $c$ and a
point $u$ in the fibre $P_{\tilde c(a)}$. The curve $d(t)$ is
called the {\em lift of the $\ro$-admissible curve $c$ through $u$
with respect to $h$} and we write from now on $c^h_u(t)$ to denote
$d(t)$. Similar to standard connection theory, we call the map
$c^h: \pi^{-1}(\tilde c(a)) \to \pi^{-1}(\tilde c(b)) : u \mapsto
c^h_u(b)$, {\em the $h$-displacement along $c$}. It is easily seen
that $c^h$ commutes with $R_g$ for $g \in G$ arbitrary, i.e.
$c^h(ug) = c^h(u)g$. Therefore, $c^h$ determines a morphism on the
fibres of $P$. The lift of a composition of $\ro$-admissible
curves, in the sense of Section \ref{sectionanchor}, equals the
composition of the corresponding $h$-admissible curves. Following
the constructions described in the previous section, we can also
consider the inverse anchored bundles of $(\nu,\ro)$ and $(\tilde
\pi_1,h)$. We have that $(c^*)^{-h} = (c^h)^{-1}$, i.e. $c^h$ is
invertible, that any $\pm h$-admissible curve projects onto a
$\pm\ro$-admissible curve and that any $\rm\ro$-admissible curve
is the projection of a $\pm h$-admissible curve. Hence, using
Theorem \ref{thmkarleafs}, we obtain $\pi(H(u))=L_{\pi(u)}$. This
result is of great importance for the development of a notion of
leafwise holonomy for principal $\ro$-lifts.

\begin{defn}\label{defholonomy}
The set of all $g \in G$ such that $ug\in H(u)$, which is called
the {\em holonomy group with reference point $u$}, is denoted by
$\Phi(u)$.
\end{defn}

The fact that $\Phi(u)$ is a subgroup follows from the following
lemma. First, note that, given any $g \in \Phi(u)$, there exists a
$\pm h$-admissible curve taking $u$ to $v=ug$, since $v \in H(u)$.
This $\pm h$-admissible curve projects onto a $\pm\ro$-admissible
loop with base point $\pi(u)=\pi(v)=x$. This implies that $v$ can
be reached from $u$ by composing $h$-admissible curves and
$(-h)$-admissible curves. Since a $h$-admissible curve is a lift
of a $\ro$-admissible curve and since a $(-h)$-admissible curve is
a lift of a $(-\ro)$-admissible curve, we obtain that $g$ is
determined by composing a finite number of $h$-displacements along
$\ro$-admissible curves and $(-h)$-displacements along
$(-\ro)$-admissible curves. In particular, using the notations
from Section \ref{sectionanchor}, we can define a map from the
loop space $C(x,N)$ to $\Phi(u)$, which is onto. These
observations are used in the proof of the following lemma.

\begin{lem}\label{lemholissubgroup}
$\Phi(u)$ is a subgroup of $G$.
\end{lem}
\begin{pf}
Given any two elements $g,g' \in \Phi(u)$ and let $ug
=((c^\ell)^{\pm h} \circ \ldots \circ (c^1)^{\pm h})(u)$, and $ug'
= ((c^{\ell+\ell'})^{\pm h} \circ \ldots \circ (c^{\ell+ 1})^{\pm
h})(u)$, for some $\pm\ro$-admissible curves $c^i$,
$i=1,\ldots,\ell+\ell'$, and where $(c^i)^{\pm h}$ stands for
$(c^i)^h$ if $c^i$ is $\ro$-admissible, and $(c^i)^{-h}$ if $c^i$
is $(-\ro)$-admissible.

Then $g'g^{-1} \in \Phi(u)$ since
\[
u{g^{-1}g'}=((c^{\ell+\ell'})^{\pm h} \circ \ldots \circ (c^{\ell+
1})^{\pm h}\circ ((c^*)^1)^{\pm h} \circ \ldots \circ
((c^*)^\ell)^{\pm h}(u),
\]
and, hence, $ug'g^{-1}$ belongs to $H(u)$. \qed
\end{pf}
In the above proof, we used the fact that any $\pm\ro$-admissible
loop $c=c^\ell\cdot \ldots \cdot c^1$ with base point $x \in M$,
we can associate a map on the fibre $\pi^{-1}(x)$ which commutes
with the right action (i.e. such a map is called an automorphism
on $\pi^{-1}(x)$). Indeed, for $u \in \pi^{-1}(x)$ and $g \in G$
arbitrary, we have
\[(c^\ell)^{\pm h} \circ \ldots \circ (c^1)^{\pm h}(ug) =
(c^\ell)^{\pm h} \circ \ldots \circ (c^1)^{\pm h}(u)g.\] Using
similar arguments as in the above proof, the set of all such
automorphisms on the fibre $\pi^{-1}(x)$ forms a group, which is
called the {\em holonomy group with reference point $x$} and
denoted by $\Phi(x)$. We thus have the following commutative
diagram:
\begin{center}
\setlength{\unitlength}{1cm}
\begin{picture}(3.84,2.50) \thicklines
\put(0.121,0.297){\rmfamily $\Phi(x)$} \put(2.63,0.297){\rmfamily
$\Phi(u)$} \put(1.12,1.95){\rmfamily $C(x,N)$}
\put(1.29,1.70){\vector(-2,-3){0.601}}
\put(1.1,0.466){\vector(1,0){1.3}}
\put(2.16,1.72){\vector(1,-2){0.458}}
\end{picture}
\end{center}

\begin{rem}\textnormal{
In the specific case where $h$ is a principal $\ro$-connection,
the situation becomes more simple. In order to define the concept
of holonomy groups it is sufficient to consider only
$\ro$-admissible loops. Indeed, if $c$ is $(-\ro)$-admissible,
then $-c$ is $\ro$-admissible, and $c^{-h} = (-c)^h$. Moreover, we
can consider reparameterisations of $\ro$-admissible curves and
the notion of $h$-displacement does not depend on the
parametrisation of $c$, in the following sense. Assume that $\phi:
[a,b] \to [c,d]$ is a diffeomorphism with $\phi(a)=c$ and
$\phi(b)=d$, then the curve $c': [c,d] \to N$, defined by
\[
c'(s) = \frac{d\phi^{-1}}{ds}(s) c(\phi^{-1}(s)),
\]
is $\ro$-admissible and, as can be seen from elementary
calculations, it follows that $h$-displacement along $c$ or $c'$
is the same. Recall the definition of the inverse $c^{-1} = -c^*$
of a $\ro$-admissible curve $c$. The following identity holds
$(c^{-1})^h = (c^h)^{-1}$.}
\end{rem}

The following properties are well known from the standard theory
of holonomy.

\begin{prop}
$(i)$ Given any $v \in H(u)$, then $\Phi(u)=\Phi(v)$. $(ii)$ Given
any $g \in G$, then $\Phi(ug)=I_{g^{-1}}(\Phi(u))$, where, $I$
denotes the inner automorphism on $G$ (i.e. for $h\in G$,
$I_{h}:G\to G:h'\mapsto hh'h^{-1}$).
\end{prop}
\begin{pf}
By definition of $H(u)$, we have that $H(ug)=R_g(H(u))$. Indeed,
$H(u)$ is the leaf of a foliation of a distribution generated by
right invariant vector fields. Thus, if $h \in\Phi(u)$, then
$h^{-1} \in \Phi(u)$ and $uh^{-1}\in H(u)$, or $H(uh^{-1})=H(u) =
H(v)$. Acting on the right by $h$, we obtain $H(u)=H(vh)$. And
since $H(u)=H(v)$, we have $h \in \Phi(v)$, proving $(i)$. Since
$H(ug)=H(u)g$, we have that, for any $h \in \Phi(u)$, then $H(uhg)
= H(ug)$. Thus $g^{-1}hg\in \Phi(ug)$, proving $(ii)$. \qed
\end{pf}

\section{Mappings between generalised connections}\label{sectiemapconnecties}
We first fix some notations. Let $(\nu',\ro')$ and $(\nu,\ro)$
denote anchored bundles with base manifolds, respectively, $M'$
and $M$ and consider an anchored bundle morphism $f:N'\to N$
between $(\nu',\ro')$ and $(\nu,\ro)$, which is fibred over $\ovl
f: M'\to M$. Assume that $\pi':P'\to M'$ and $\pi: P \to M$ are
principal fibre bundles with structure groups, respectively $G'$
and $G$. Furthermore, we assume that a principal fibre bundle
morphism $F: P' \to P$ between $P'$ and $P$ is given, such that
$F$ is also fibred over the map $\ovl f: M'\to M$ between the base
spaces. The group morphism between $G'$ and $G$, corresponding to
$F$, is denoted by $\ovl F: G'\to G$, i.e. for all $u' \in P'$ $g'
\in G'$, we have $F(u'g')= F(u')\ovl F(g')$.

The principal fibre bundle morphism $F$ is called a {\em morphism
between the principal $\ro'$-lift $h'$ and the principal
$\ro$-lift $h$} if the map $(F,f)$, defined by $(F,f):(\pi')^*N'
\to \pi^*N:(u',s')\mapsto (F(u'),f(s')$, is an anchored bundle
morphism between $(\tilde \pi'_1,h')$ and $(\tilde \pi_1,h)$. More
precisely we have that:
\[ TF \left(h'(u',s') \right) = h(F(u'),f(s')).
\]

\begin{thm} \label{thmreductionoflifts}
Assume that $f$ is an isomorphism, and that $F$ is a principal
fibre bundle morphism from $P'$ to $P$, fibred over $\ovl f$. Let
$h'$ be a principal $\ro'$-lift on $P'$. There exists a unique
principal $\ro$-lift $h$ such that $F$ is a morphism between $h'$
and $h$. The holonomy group $\Phi(u')$ of $h'$ is mapped by $\ovl
F$ onto $\Phi(F(u'))$.
\end{thm}
\begin{pf}
Let $u$ denote an arbitrary point of $P$, with $\pi(u)=x$. Then
fix an element $u'$ in $P'_{\ovl f^{-1}(x)}$ and an element $g$ in
$G$ such that $F(u')=ug$. Define $h(u,s) \in T_uP$, for any $s \in
N_{\ovl f^{-1}(x)}$, by
\[
h(u,s) = TR_{g^{-1}} \left(T_{u'} F (h'(u',f^{-1}(s)))\right).
\]
This tangent vector in $T_uP$ is well defined, in the sense that
it does not depend on the choice of $u'$, since for any other
element $v'=u'g'$, then $v'$ satisfies $F(v')=F(u')\ovl F(g')= uh$
with $h=g\ovl F(g')$, implying that
\begin{eqnarray*}
h(u,s) &=& TR_{h^{-1}} \left(T_{v'} F(h'(v',f^{-1}(s)))\right)\\
&=& T R_{h^{-1}} \left(T_{v'} F(TR_{g'}h'(u',f^{-1}(s)))\right)\\
&=& TR_{g^{-1}}TR_{\ovl F({g'}^{-1})} TR_{\ovl F(g')} \left(T_{u'}
F(u',f^{-1}(s)) \right) \\
&=&TR_{g^{-1}} \left(T_{u'} F (h'(u',f^{-1}(s)))\right).
\end{eqnarray*}
In this way, we have constructed a mapping $h:\pi^*P \to TP$,
which is clearly right invariant and, by definition, it follows
that $F$ is an anchored bundle morphism between
$(\tilde\pi'_1,h')$ and $(\tilde \pi_1,h)$. From the fact that
$f^{-1}$ maps any $\pm\ro$-admissible curve onto a
$\pm\ro'$-admissible curve, we have that $H(u')$ is mapped by $F$
onto $H(F(u'))$, concluding the proof. \qed
\end{pf}

In the specific case that $P'$ is a reduced subbundle of $P$, i.e.
$F$ is an injective immersion and $\ovl F$ is an monomorphism,
then we say that $h$ is {\em reducible to a principal $\ro$-lift
on $P'$}. This is important for our treatment of holonomy, where
we prove a generalisation of the Reduction Theorem, which says
that $H(u)$ is a reduced subbundle with structure group the
holonomy group $\Phi(u)$ and that $h$ is reducible to $H(u)$.

For the following theorem we take for $(\nu',\ro')$ the pull-back
anchored bundle of $(\nu,\ro)$ under $i:L_x\hookrightarrow M$,
with $L_x$ the leaf through some $x \in M$. Let $P'=i^*P$ and $F:
P' \to P$ the projection onto the second factor. Note that the
structure group of $P'$ is precisely $G$.

\begin{thm} \label{thmpullbackcon}
There exists a unique $\ro'$-lift $h'$ on $P'$ such that $F$ is a
morphism between $h'$ and $h$. Moreover, $F(H(u')) = H(F(u'))$
and, therefore $\Phi(u')=\Phi(F(u'))$.
\end{thm}
\begin{pf}
Since $F$ is an injective immersion, we know from Section
\ref{sectionanchor} that a unique anchor map $h'$ on $P'$ can be
defined such that $F$ is an anchored bundle morphism between
$(\tilde \pi'_1,h')$ and $(\tilde \pi_1,h)$. It is trivial to
check that $h'$ satisfies the ``right invariance'' condition
making it into a principal $\ro$-lift.

The fact that the induced foliations coincide follows from the
fact that $\pm\ro$-admissible curves are in one-to-one
correspondence with the $\pm\ro'$-admissible curves.  \qed
\end{pf}
In the following section we prove that the holonomy groups
$\Phi(u)$ of a principal $\ro$-lift is a Lie subgroup of $G$. In
view of the above theorem, we will assume that, without loss of
generality, we are working with the $\ro'$-lift $h'$ on the bundle
$i^*P$, with $i:L_x\hookrightarrow M$. Indeed, the holonomy groups
of $h$ and $h'$ are the same.

\section{Leafwise Holonomy of a principal $\ro$-lift}\label{sectieholonomy}

In view of the above comment, we have that $M=L_x$ is a connected
manifold and that $\tilde D = TM$. The main consequence of these
assumptions is that the distribution $\tilde Q$ generated by a
principal $\ro$-lift $h$ is regular, i.e. has constant rank. We
have to prove that $\dim \tilde Q_u= \dim \tilde Q_v$, given two
arbitrary points $u,v$ in $P$. Let $x=\pi(u)$ and $y=\pi(v)$.
Then, since $M=L_x$, there exists a composite flow $\Phi$
associated with $(\ro\circ\s^\ell,\ldots,\ro\circ\s^1)$ of vector
fields in $\cal D$ and a composite flow parameter $T$ such that
$\Phi_T(x)=y$ (cf. Theorem \ref{thmsus}). Consider the vector
fields $(\s^i)^h$ in $\cal Q$. The flows of $(\s^i)^h$ and $\ro
\circ \s^i$ are $\pi$-related by definition, and therefore, if
$\Phi^h$ is the composite flow of $((\s^\ell)^h,\ldots,(\s^1)^h)$,
we have $\pi(\Phi^h_T(u))=y$, or there exists a $g \in G$ such
that $\Phi^h_T(u)g=v$. By definition of $\tilde Q$ we have
$T\Phi^h_T(\tilde Q_u)=\tilde Q_{\Phi^h_T(u)}$. On the other hand
since $\cal{D}^h$ consists of right invariant vector fields and
since $\cal D^h$ generates $\tilde Q$, we have $TR_h(\tilde
Q_w)=\tilde Q_{wh}$ for any $w \in P$ and $h\in G$. Thus, we
obtain $TR_g\circ T\Phi^h_T$ is an isomorphism from $\tilde Q_u$
to $\tilde Q_v$.

Take an arbitrary point $u\in P$ and consider the linear subspace
$\alg{g}(u)$ of $\alg{g}$ defined by $T_e\s_u(\alg{g}(u)) = V_u\pi
\cap \tilde Q_u$.

\begin{prop}
Take $u\in P$ and let $g\in G$ arbitrary. Then, $(i)$ $\alg{g}(u)
= \alg{g}(v)$ for any $v \in H(u)$, $(ii)$ $Ad_{g^{-1}}(
\alg{g}(u)) = \alg{g}(ug)$ and $(iii)$ $\alg{g}(u)$ is a Lie
subalgebra of $\alg{g}$.
\end{prop}
\begin{pf}
$(i)$ follows from the fact that $V\pi$ and $\tilde Q$ are
invariant under the image of the tangent map to a composite flow
associated with vector fields in $\cal{D}^h$. $(ii)$ follows from
$TR_g \circ T_e\s_u = T_e\s_{ug} \circ Ad_{g^{-1}}$, $TR_g
(V_u\pi) = V_{ug}\pi$ and $TR_g(\tilde Q_u) = \tilde Q_{ug}$.
$(iii)$ follows from $[\s(A),\s(B)]=\s([A,B])$, for $A,B\in
\alg{g}$ and the fact that $\tilde Q$ is involutive (since it is
integrable, by definition). \qed
\end{pf}

These properties allow us to consider the connected Lie group
$\Phi^0(u)$ generated by the Lie algebra $\alg{g}(u)$, which is
called the {\em restricted holonomy group}. From the preceding
proposition, we have that $\Phi^0(u)=\Phi^0(v)$ for $v \in H(u)$
and $\Phi^0(ug)=I_{g^{-1}}(\Phi^0(u))$.

We prove that $\Phi^0(u)$ is a normal subgroup of $\Phi(u)$ and
that $\Phi(u)/\Phi^0(u)$ is countable, implying that $\Phi(u)$ is
a Lie-subgroup of $G$ whose identity component is precisely
$\Phi^0(u)$, see \cite[p 73]{koba}. We first prove that
$\Phi^0(u)$ is normal subgroup of $\Phi(u)$.

Let $h \in \Phi^0(u)$. By construction of the Lie subgroup
$\Phi^0(u)$ (i.e. it is the leaf through $e$ of the left invariant
distribution generated by $\alg{g}(u)$), $h$ is obtained from $e$
by a composite flow associated with left invariant vector fields
generated by $\alg{g}(u)$. Note that, if $g(t)$ denotes the
integral curve through $e$ of the left invariant vector field
corresponding to $A \in \alg{g}(u)$, then $ug(t) \in H(u)$, since
$\s(A)$ determines a vector field tangent to $H(u)$, and hence
$g(t) \in \Phi(u)$. We therefore have $\Phi^0(u) < \Phi(u)$. Since
$\Phi^0(ug) = I_{g^{-1}}(\Phi^0(u))$ and $\Phi^0(u) = \Phi^0(ug)$
for any $g \in \Phi(u)$ (i.e. $\alg{g}(u) = \alg{g}(ug)$), we may
conclude that $\Phi^0(u)$ is a normal subgroup of $\Phi(u)$.

Following a similar reasoning as in \cite[p 73]{koba}, we now
prove that $\Phi(u)/ \Phi^0(u)$ is countable by constructing a
group morphism from $\pi_1^N(L_x)$ to $\Phi(u)/\Phi^0(u)$ which is
onto. Since $\pi_1^N(L_x) < \pi_1(M)$ and $\pi_1(M)$ is at most
countable, the obtain that the quotient is also countable.

\begin{pf}
Let us first make the following basic observation. In order to
prove that the map between $C(x,N)$ and $\Phi(u)$ reduces to a
well defined morphism from $\pi_1^N(L_x) \to \Phi(u)/\Phi^0(u)$,
we must prove that the images of two $\pm\ro$-admissible loops,
whose base curves are homotopic, equal up to an element in
$\Phi^0(u)$. This is achieved by using some results from standard
connection theory. Once we have obtained this morphism
$\pi_1^N(L_x) \to \Phi(u)/\Phi^0(u)$ it is easily seen to be onto,
which concludes the proof.

Consider a connection $\w$ on $P$, such that $\im h^\w$ is a
subspace of $\tilde Q$. This is always possible since $\tilde Q$
is regular and $T\pi(\tilde Q)=TM$. Consider the coefficient
$\chi$ of $h$ with respect to $\w$ (see Section
\ref{sectieprincipalcon}). Note that $T_e\s_u(\chi(u,s)) = h(u,s)
- h^\w(u,\ro(s))$ is contained in $\tilde Q$ for any $(u,s)\in
\pi^*N$. This implies that $\chi(u,s) \in \alg{g}(u)$, for all $s
\in N_{\pi(u)}$. On the other hand, the holonomy group with
reference point through $x$ of the standard connection $\w$ is a
subgroup of $\Phi(u)$ and the restricted holonomy group of $\w$ is
a subgroup of $\Phi^0(u)$, since the smallest integrable
distribution spanned by $\im h^\w$ must be contained in $\tilde Q$
(see \cite{koba}).

In Section \ref{sectieprincipalcon} we have proven that the
$h$-lift $c^h_u(t)$ of a $\ro$-admissible curve through $u \in
\pi^{-1}(x)$ equals $c^h_u(t) = d^{\w}(t) g(t)$, where $g(t)$ is a
curve in $G$ with $g(a)=e$ and $R_{g(t)^{-1}} \dot g(t) =
\chi(d^\w(t),c(t))$, and where $\dot d^\w(t)
=h^\w(d^\w(t),\dot{\tilde c}(t))$, with $d^\w(a)=u$. In particular
we have $g(t) \in \Phi^0(u)$ (since the image of $\chi$ is
contained in $\alg{g}(u)$). This is also valid for the inverted
anchored bundles. Thus we can conclude that any element belonging
to $\Phi(u)$ can be written as a product of elements belonging to
the holonomy group of $\w$ at $u$ and of elements in $\Phi^0(u)$.
Moreover, if the base of a $\pm\ro$-admissible curve is homotopic
to zero, then the corresponding product of elements is entirely
contained in $\Phi^0(u)$, since the restricted holonomy group of
$\w$ is a subgroup of $\Phi^0(u)$. This completes the proof. \qed
\end{pf}

\begin{cor}
The holonomy group $\Phi(u)$ is a Lie subgroup of the structure
group $G$ with Lie algebra $\alg{g}(u)$.
\end{cor}

We are now able to state a generalisation of the reduction theorem
for principal $h$-lifts.
\begin{thm}
$H(u)$ is a reduced subbundle of $P$ with structure group
$\Phi(u)$ and $h$ reduces to a principal $\ro$-lift on $H(u)$.
\end{thm}
\begin{pf}
It is sufficient to prove that, given a point $y \in L_x$, there
exists a neighbourhood $U \ni y$ and a section $\s$ of $P$ defined
on $U$ such that $\s(U) \subset H(u)$. The existence of such a
cross-section follows by using a result from \cite[p 84]{koba}
with respect to a connection $\w$ with horizontal distribution
contained in $\tilde Q$.

Since $H(u)$ is the leaf of the foliation induced by ${\cal Q}$,
we can consider the pull-back anchor map of $h$. Using the fact
that $H(u)$ is a principal fibre bundle over $L_{\pi(u)}$ and
using Theorem \ref{thmreductionoflifts}, it is easily seen that
$h$ is reducible to the pull-back of $h$. \qed
\end{pf}

Assume that $\dim M \ge 2$. Then, since $H(u)$ is connected, there
exists a standard principal connection $\ovl \w$ on $H(u)$ whose
holonomy group is the structure group $\Phi(u)$ (see \cite[p
90]{koba}). Using Theorem \ref{thmreductionoflifts} from Section
\ref{sectiemapconnecties}, then $\w$ can be extended to a
connection on $P$.
\begin{cor} If $\dim M\ge 2$, then there exists a connection $\w$ on $P$ such that the
holonomy groups of $\w$ equal the holonomy groups of the lift $h$.
\end{cor}

\section{Possible field of applications}

The equations of motion a free particle subjected to linear
nonholonomic constraints can be described as the ``geodesics'' of
a unique connection along the natural injection of the constraint
distribution into the tangent bundle of the configuration
manifold, see \cite{mijzelf}. This unique generalised connection
admits a notion of holonomy and, consequently, one can wonder
wether the holonomy groups may play a role in the study of
nonholonomic motions.

Another field of application could be found in sub-Riemannian
geometry, see \cite{sub}. However, until now, we haven't been able
to construct a unique generalised connection in sub-Riemannian
geometry. These possible applications of the above developed
theory on holonomy groups of generalised connections is left for
future work.

\section*{Acknowledgements} This work has been supported by a grant from the
``Bijzonder Onderzoeksfonds'' of Ghent University. I am indebted
to F. Cantrijn for the many discussions and the careful reading of
this paper.

%GATHER{C:\program files\miktex\bibtex\bib\bibliotheek\Bibliotheek.bib}   % For Gather Purpose Only
%GATHER{holonomy2earchive.bbl} % For Gather Purpose Only
\providecommand{\bysame}{\leavevmode\hbox
to3em{\hrulefill}\thinspace}
\providecommand{\MR}{\relax\ifhmode\unskip\space\fi MR }
\providecommand{\MRhref}[2]{%
  \href{http://www.ams.org/mathscinet-getitem?mr=#1}{#2}
} \providecommand{\href}[2]{#2}

\end{document}